\documentclass [11pt]{article}
\usepackage{amssymb,amsmath,graphicx}
\newcommand{\eps}{\varepsilon}
\newcommand{\RR}{\mathbb{R}}

\newcommand{\ZZ}{\mathbb{Z}}

\newcommand{\f}{\frac}
\newcommand{\Res}{\mathrm{Res}}
\renewcommand{\Re}{\mathrm{Re}}
\renewcommand{\Im}{\mathrm{Im}}
\newcommand{\pd}{\partial}

\newcounter{problem}

\newcommand{\dst}{\displaystyle}
\newcommand{\unl}{\underline}
\newtheorem{lem}{Lemma}
\newcommand{\pf}{\noindent{\it Proof}. }
\newcommand{\eop}{\hfill $\Box\quad$}

\newcommand{\ba}[1]{\begin{array}{#1}}
\newcommand{\ea}{\end{array}}
\newcommand{\beq}[1]{\begin{equation}\label{#1}}
\newcommand{\eeq}{\end{equation}}
\newcommand{\num}[1]{{\rm(\ref{#1})}}

\title{Asymptotics of alternating harmonic series with attenuation}
\author{Sergey Sadov\footnote{
E-mail: serge.sadov@gmail.com}}
\date{}

\begin{document}
\maketitle

\begin{abstract}
%
%\begin{verse}{\small}
%\noindent
%{\bf Abstract
We find the asymptotics of the series $\sum_{n=1}^\infty (-1)^n n^{-1} \exp(-t/n)$ as $t\to+\infty$. 
The answer is an oscillating function of $t$ dominated by $\exp(-(2\pi t)^{1/2})$.
The intermediate step is to find the
asymptotics of the two-dimensional Fourier transform $\hat F(\xi)$ of the function $F(x)=(1+\exp(\|x\|^2))^{-1}$ as $\|\xi\|\to\infty$. 


\end{abstract}
%\end{verse}

%\begin{center}{
%MATH 6101 --- Asymptotic Analysis --- Fall 2011
%}
%\end{center}

%\subsection*{Solution of Problem 7 of Asst.4} 

The problem discussed in this note is an yet another example of a challenge born of a teaching mishap. 
I offered it by mistake among a set of exercises on the Euler-Maclaurin formula in a course of asymptotic analysis 
at Memorial University of Newfoundland in the Fall 2011. 

%I am not aware of a reference that would treat asymptotic problems of such a type.  

\bigskip\noindent
{\bf Problem.}\@
Find asymptotics of the series
\beq{S}
 S(t)=\sum_{n=1}^\infty \f{(-1)^n}{n}\,e^{-t/n}
\eeq
as $t\to+\infty$.

\bigskip
In the interests of those who want to take up the challenge the answer is only given at the end of the note.  
%(p.~\pageref{page:answer})
The solution below is long but detailed. It is intended to be understood by an asymptotic analysys course student who wants to 
get their hand dirty with the Saddle Point Method.

\bigskip\bigskip\noindent
{\bf Solution.}\@

\medskip
\noindent\unl{Part 1: Derivation of integral representation}

\medskip
Let us consider a more general series
\beq{Sgen}
 S(z,\nu,t)=\sum_{n=1}^\infty \f{z^n}{n^{\nu}}\,e^{-t/n}.
\eeq
%We will always aasume that $t\geq 0$.

\begin{lem}\label{lem:conv}
(a) The series \num{Sgen} converges absolutely if $|z|<1$ for any $\nu\in\RR$. 
\\[0.75ex]
(b) If $\nu>1$, then the series converges absolutely and uniformly in the closed unit circle $|z|\leq 1$. 
\\[0.75ex]
(c) 
If $0<\nu\leq 1$, $|z|=1$, $z\neq 1$,
the series converges conditionally. Also, for any $\eps\in(0,2)$, the convergence is uniform in the region
$$
D_\eps=\left\{z\,:\, |z|\leq 1,\; \Re\, z\leq 1-\eps\right\}.
$$
%(the unit disk with deleted segment ajacent to the point $z=1$). 
%\\[0.5ex]
(d) Consequently, if $|z|=1$, $z\neq 1$, then
\beq{Sradlim}
S(z,\nu,t)=\lim_{\rho\to 1^-}S(z\rho,\nu,t).
\eeq
\end{lem}

\pf
Parts (a) and (b) are obvious since the series is dominated by
$\sum |z|^n n^{-\nu}$ in the case (a) and by the $z$-independent sum $\sum n^{-\nu}$
in the case (b).

\smallskip
To prove (c), let us write $e^{-t/n}=1-\delta_n$,
where $\delta_n=O(n^{-1})$ (for a fixed $t$).
The series $\sum z^n n^{-\nu} \delta_n$ convereges absolutely
and uniformly in $\{|z|\leq 1\}$ (similarly to the series 
$S(z,\nu+1,t)$). 
It remains to establish the convergence properties (c)
for the series
\footnote{The function $\mathrm{Li}_\nu(z)$ defined by this series is
called {\em polylogarithm of order $\nu$}.} 
$$
 \sum_{n=1}^\infty n^{-\nu} z^n.
$$ 
It is a standard application of Dirichlet's test:
the sequence $\{n^{-\nu}\}$ decreases and the uniform in $D_\eps$ bound for the partial sums 
$$
 \left|\sum_{n=1}^N z^n\right|\leq\f{2}{|1-z|}\leq\f{2}{\eps}
$$ 
holds. 

\smallskip
To prove (d), write $z=e^{i\theta}$,where $\theta\in(0,2\pi)$. Let $\eps>0$ be such that
$\eps<1-\cos\theta$. Then $\rho z\in D_\eps$ for any $\rho\in[0,1]$. The series 
$$
 \sum_{n=1}^\infty \f{(\rho e^{i\theta})^n}{n^\nu} e^{-t/n}
$$
converges uniformly w.r.to $\rho$. Hence termwise passing to the limit as $\rho\to 1^{-}$ is justified, and \num{Sradlim} follows.
\eop

\iffalse
\medskip\noindent
{\bf Remark.}
In the special case $z=-1$, the properties (c)
are more easily seen by the alternating series test. 
The sequence $n^{-\nu}e^{-t/n}$ is decreasing for 
$n>t/\nu$: indeed, setting $t/n=h$, we have
$$
 \f{d}{dh}\log(h^\nu e^{-h})=\f{\nu}{h}-1>0
$$ 
when $h<\nu$.

Clearly, for any $r\in[0,1]$ the sequence $r^n n^{-\nu}e^{-t/n}$ is also decreasing when $n>t/\nu$. 
Therefore, if $N>t\nu$, we have
$$
\left|\sum_{n=N}^\infty \f{(-r)^n}{n^\nu}e^{-t/n}\right|
<\f{1}{N^\nu},
$$
which means that the series with parameter $r\in[0,1]$ converges uniformly. In particular,
$$
 S(t)=\lim_{r\to 1^{-}} \sum_{n=1}^\infty \f{(-r)^n}{n} e^{-t/n}.
$$
\fi

\bigskip\noindent
{\bf Remark}. In Lemma~\ref{lem:conv}, the parameter $t$
can be any complex number. For fixed $z$ ($|z|\leq 1,\; z\neq 1$) and $\nu>0$, the function $S(z,\nu,t)$ is an entire analytic function of $t$. 

It is interesting also to note the relations
$$
 \ba{l}
  \dst \f{\pd S(z,\nu,t)}{\pd t}=-S(z,\nu+1,t),
\\[2ex]\dst
  \f{\pd S(z,\nu+1,t)}{\pd z}=z^{-1} S(z,\nu,t),
\\[2ex]\dst
  \f{\pd^2 S(z,\nu,t)}{\pd t\,\pd z}=-z^{-1} S(z,\nu,t).
 \ea
$$
\eop 

\bigskip\noindent
Recall the series representing Bessel function of order $n\geq 0$:
\beq{bes}
 J_n(u)=\sum_{k=0}^\infty \f{(-1)^k (u/2)^{2k+n}}{k!\,\Gamma(k+n+1)}.
\eeq

\begin{lem} For $\nu\geq 1$ and $|z|<1$, the function $S(z,\nu,t)$ has the integral representation
\beq{irepz}
S(z,\nu,t)=\int_0^\infty \f{ze^{-x}}{1-ze^{-x}}\,\left(\f{x}{t}\right)^{\f{\nu-1}{2}}\,J_{\nu-1}(2\sqrt{tx})\,dx.
\eeq
The same representation remains valid if $|z|=1$,
$z\neq 1$.
In particular,
\beq{irep1}
S(t)=S(-1,1,t)=-\int_0^\infty \f{1}{e^x+1}\,J_0(2\sqrt{tx})\,dx.
\eeq

%where 
%$$
% J_0(u)=\sum_{n=0}^\infty \f{(-u^2/4)^n}{(n!)^2}
%$$
%is the Bessel function of order 0.
\end{lem}

\pf
Suppose $|z|<1$. 
Since 
$$
 \f{\Gamma(k+\nu)}{n^{k+\nu}}=\int_0^\infty e^{-nx}\,x^{k+\nu-1}\,dx,
$$
we have
$$
\ba{rcl}
 S(z,\nu,t)&=&\dst
\sum_{n=1}^\infty z^n\sum_{k=0}^\infty \f{(-t)^k}{k!}\,n^{-\nu-k}
\\[3ex]&=&\dst
\sum_{n=1}^\infty z^n\sum_{k=0}^\infty \f{(-t)^k }{k!\,\Gamma(k+\nu)}\,\int_0^\infty e^{-nx}\,x^{k+\nu-1}\,dx
\\[3ex]&=&\dst
\sum_{k=0}^\infty \f{(-t)^k }{k!\,\Gamma(k+\nu)}\,\int_0^\infty\left(\sum_{n=1}^\infty z^n  e^{-nx}\right)\,x^{k+\nu-1}\,dx
\\[3ex]&=&\dst
\sum_{k=0}^\infty \f{(-t)^k}{k!\,\Gamma(k+\nu)}\,\int_0^\infty\f{ze^{-x}}{1-ze^{-x}}\,x^{k+\nu-1}\,dx
\\[3ex]&=&\dst
\int_0^\infty\f{ze^{-x}}{1-ze^{-x}}\,x^{\nu-1}\,
\sum_{k=0}^\infty \f{(-tx)^k}{k!\,\Gamma(k+\nu)}\,dx.
\ea
$$
In the above transformations, convergence is absolute throughout,
so we can interchange the order of summations and integration freely.

Setting $u=2\sqrt{tx}$ and
comparing the integrand in the RHS with \num{bes},
we see that
$$
 x^{\nu-1}\,
\sum_{k=0}^\infty \f{(-tx)^k}{k!\,\Gamma(k+\nu)}
=\left(\f{x}{t}\right)^{\f{\nu-1}{2}}\,\sum_{k=0}^\infty \f{(-1)^k (u/2)^{2k+{\nu-1}}}{k!\,\Gamma(k+\nu)}
=\left(\f{x}{t}\right)^{\f{\nu-1}{2}}\,J_{\nu-1}(u),
$$
so that \num{irepz} follows.

Denote temporarily the RHS of \num{irepz} by $S_*(z,\nu,t)$.
The integral in \num{irepz} converges absolutely if
$|z|\leq 1$, $z\neq 1$. (Convergence at $x\to\infty$ takes
place for any $z$ and convergence at $x\to 0^+$ is certain
if $\nu\geq 1$, $z\notin [1,\infty)$.) 
Consequently, if $|z|=1$, $z\neq 1$,
$$
 \lim_{\rho\to 1^-} S_*(z\rho,\nu,t)=S_*(z,\nu,t).
$$
On the other hand, Lemma~\ref{lem:conv} asserts the same
for $S(z,\nu,t)$. Thus, \num{irepz} remains valid
in the case  $|z|=1$, $z\neq 1$ by continuity.
\eop

\bigskip\noindent
{\bf Corollary}.
Let 
\beq{tlam}
 \lambda=2\sqrt{t}
\quad\mbox{\rm and}\quad
S(t)=S_*(\lambda).
\eeq
Taking $\sqrt{x}$ as the new variable of integration in
\num{irep1}, we get
\footnote{$S_*(\lambda)$ is called the {\it Hankel transform}\ of the function $-2x/(\exp(x^2)+1)$.}
\beq{ireplam}
 S_*(\lambda)=-\int_0^\infty J_0(\lambda x)\f{2x}{e^{x^2}+1}\,dx.
\eeq

\bigskip
\noindent\unl{Part 2: Transformations of the integral \num{ireplam}}

\medskip
Substituting the integral representation
$$
 J_0(u)=\f{1}{2\pi}\int_0^{2\pi} e^{iu\cos\theta}\,d\theta, 
$$
into \num{ireplam} we express $S_*(\lambda)$ as a double integral
\beq{dblint}
S_*(\lambda)=-\f{1}{2\pi}\int_0^\infty \int_0^{2\pi}
e^{i\lambda x\cos\theta}\,d\theta\, \f{2x}{e^{x^2}+1}\,dx.
\eeq
Experimental study suggests that $S_*(\lambda)=O(e^{-\lambda})$. Unfortunately, it is not clear how a precise asymptotics can be derived from the representation \num{dblint}, because the edge $x=0$ interferes in all attempts: integration by parts or
rotation of the integration path in the complex $x$-plane.
We will not describe such attempts in detail, since
the successful solution will be based on a different transformation of the integral \num{ireplam}. Note only
that if instead of the factor $2x$ we had an even function of
$x$ in the integrand of \num{dblint}, then it would be possible to extend domain of integration w.r.t.\ $x$ to
$\RR$, getting rid of the edge; the saddle point method would
then be applicable in a standard manner. 

\medskip
The representation we need will come from a known general 
formula of Fourier analysis.

\begin{lem}
Suppose $\int_0^\infty |f(r)| r\,dr <\infty$.
Let $r=r(x,y)=(x^2+y^2)^{1/2}$. Then the 2-dimensional Fourier transform of $F(x,y)=f(r(x,y))$ is
\beq{FT}
 \hat F(\xi,\eta)\equiv \int_{\RR^2} F(x,y) e^{-ix\xi-i y\eta}\,dx\,dy
\;=\;
2\pi \int_0^\infty  J_0(r\rho)\,f(r)\,r\,dr,
\eeq
where $\rho=(\xi^2+\eta^2)^{1/2}$.
\end{lem}

\pf
Using the polar coordinates $(r,\theta)$ in the $(x,y)$ plane
and  $(\rho,\phi)$ in the $(\xi,\eta)$ plane, we get
$$
\hat F(\xi,\eta)=\int_{0}^\infty \int_0^{2\pi} f(r)\, e^{-ir\rho\cos(\phi-\theta)}\,d\theta\,r\,dr.
$$
Clearly, the change of variable $\theta\to \theta-\phi$
shows that the integral does not depend on $\phi$;
thus 
$$
\ba{rcl}
\hat F(\xi,\eta)&=&\dst
\int_{0}^\infty f(r)\,r\,\int_0^{2\pi} f(r)\, e^{ir\rho\theta}\,d\theta\,dr
\\[3ex]&=&\dst
\int_{0}^\infty f(r)\,r\cdot\,2\pi\,J_0(r\rho)\,dr,
\ea
$$
as stated.
\eop

\bigskip\noindent
{\bf Corollary}.
The integral \num{ireplam} can be written as the 2-dimensional Fourier integral
\beq{ifou}
S_*(\lambda)=-\f{1}{\pi}\int_{-\infty}^\infty\int_{-\infty}^\infty
\f{e^{i\lambda x}}{1+e^{x^2+y^2}}\,dx\,dy.
\eeq
Indeed, it suffices to replace $r$ in \num{ireplam} by $r$ and to use \num{FT} with $f(r)=2/(e^{r^2+1}+1)$, $(\xi,\eta)=(-\lambda,0)$.

\bigskip\noindent
{\bf Remark}. Expanding $(1+e^{x^2+y^2})^{-1}$ as a geometric
series and integrating separately w.r.t. $x$ and $y$, one
obtains an independent proof of the fact that $S(t)$ equals 
the RHS of \num{ifou} when $\lambda=2\sqrt{t}$. 

\bigskip
\noindent\unl{Part 3: Asymptotic analysis of the integral \num{ifou}}

Asymptotics of the integral \num{ifou} is determined by the
complex singularities of the integrand. 
As a preparation, let us locate complex zeros in the $z$-plane of the denominator
$$
 Q(z,y)=1+e^{z^2+y^2}.
$$
where $z=x+iu$, and $x$, $u$ ,$y$ are real.

The equation $Q(z,y)=0$ is equivalent to
$$
 z^2+y^2=i(2k+1)\pi, \qquad k\in\ZZ.
$$
Or, in the real form,
\beq{zQ}
 \ba{l}
 \dst x^2+y^2=u^2,
\\[2ex]
\dst
  x u =\pi\left(k+\f{1}{2}\right).
 \ea
\eeq

Precise asymptotic analysis of the integral \num{ifou} is significantly more complicated than singularities-based analysis of Fourier integrals in the one-dimensional case. 
As a first, rather easy step, we obtain a rough exponential $o$-estimate.  

\begin{lem}
For any $a<\sqrt{\pi/2}$ we have
\beq{roughO}
S_*(\lambda)=o(e^{-a\lambda}).
\eeq
\end{lem}

\pf
Clearly, for any fixed $y>0$ and $h>0$
$$
|1+e^{(x+ih)^2+y^2}|\to \infty 
$$
as $x\to\pm\infty$.
Therefore, assuming $Q(x+ih,y)\neq 0$, $\forall x\in\RR$, we
have
\beq{intxres}
 \int_{-\infty}^\infty
\f{e^{i\lambda x}\,dx}{Q(x,y)}
=
 \int_{-\infty}^\infty
\f{e^{i\lambda (x+ih)}\,dx}{Q(x+ih,y)}
+\,2\pi i\,\sum
\Res\,\f{e^{i\lambda z}}{Q(z,y)},
\eeq
where the sum of residues is taken over all $z=x+iu$ such that $0<u<h$ and $Q(z,y)=0$.

The equation  
$Q(x+iu,y)=0$ does not have solutions with $x,y\in\RR$,
$0\leq u<\sqrt{\pi/2}$.
Indeed, looking at the system \num{zQ}, 
we see that
$|x|\leq u$ and $|\pi/2|\geq |xu|\geq u^2$.

Therefore, for $h<\sqrt{\pi/2}$, the sum in the RHS of \num{intxres} is void. Thus
$$
S_*(\lambda)=-\f{1}{\pi}\int_{-\infty}^\infty 
\int_{-\infty}^\infty \f{e^{i\lambda (x+ih)}}{1+e^{(x+ih)^2+y^2}}\,dx\,dy= O(e^{-\lambda h}).
$$
Given $a\in(0,\sqrt{\pi/2})$, the estimate \num{roughO} follows by taking $h\in(a,\sqrt{\pi/2})$.
\eop

\medskip
We want eventually to nail down an asymptotical term
of the exponential order $O(e^{-\lambda/\sqrt{2}})$ precisely. 

The procedure consists of simple steps, but it is rather delicate overall.
We need to take into consideration some poles of 
$Q(z,y)^{-1}$ to obtain a nonempty the residue sum
in the RHS \num{intxres}.
On the other hand, we must avoid stepping on a pole 
during the double integration (in both $x$ and $y$).
Moreover, any estimates we get while keeping $y$ fixed
should be uniform or explicit enough to justify subsequent integration with respect to $y$.
 
The following lemma will allow us to control the integrand
in the double integral uniformly.

\begin{lem}\label{lem:Qlb}
Suppose $u>0$ and a closed set $K\subset\RR$ are
such that the equation $Q(x+iu,y)=0$ does not have solutions
$x\in\RR$, $y\in K$. Then there exists $\alpha=\alpha(a,K)>0$
such that 
\beq{lbQ}
|Q(x,y)|\geq \alpha e^{x^2+y^2}\qquad \forall x\in\RR,\;y\in K.
\eeq
\end{lem}

\pf 
We will estimate $|Q(x,y)|$ from below in two cases separately: first for large
$x^2+y^2$, then in a bounded region of the $(x,y)$ plane. 

\smallskip
1) Suppose that $x^2+y^2>u^2+\ln 2$. Then $2e^{u^2}<
e^{x^2+y^2}$, hence
$$
\ba{rcl}
 |Q(x,y)|&\geq &\dst
 |e^{\Re(x+iu)^2+y^2}-1|=e^{-u^2}(e^{x^2+y^2}-e^{u^2})
\\[1.5ex]
&\geq &\dst
\f{e^{-u^2}}{2}\,e^{x^2+y^2}.
\ea
$$

\smallskip
2) The set 
$$
\Omega=\{(x,y)\,|\, x^2+y^2\leq u^2+\ln 2,\; y\in K\}
$$
is a bounded closed subset of $\RR^2$, hence compact.
By assumption, $Q(x,y)\neq 0$ when $(x,y)\in\Omega$.
The function $e^{x^2+y^2}/Q(x,y)$ is continuous in $\Omega$, hence bounded: $|e^{x^2+y^2}/Q(x,y)|\leq C$.

\smallskip
The inequality \num{lbQ} with 
$$
\alpha=\min\left\{\f{e^{-u^2}}{2},\;\f{1}{C} \right\}
$$
follows.
\eop

Next comes the crusial step. We choose a splitting parameter
$b$ for integration with respect to $y$ and write \num{ifou}
in the form
\beq{isplit}
 -\pi S_*(\lambda)=\int_{|y|>b} T(y,\lambda)\,dy+
\int_{|y|\leq b} T(y,\lambda)\,dy\;=\;I_1+I_2,
\eeq
where
\beq{intt}
T(y,\lambda)=\int_{-\infty}^\infty
\f{e^{i\lambda x}}{1+e^{x^2+y^2}}\,dx.
\eeq

We are going to make a complex shift $x\to x+ia$ in the
integral \num{intt}. However, we will use two different
values of $a$ in the two cases, $|y|\geq b$ and $|y|\leq b$.

In the case $|y|\geq b$, we will take $a=a_1$ so as to 
ensure that no zeros of $Q(z,y)$ lie in the strip
$0\leq \Im z\leq a_1$. 
It will be possible to choose
$a_1>\sqrt{\pi/2}$ and to obtain the estimate
$I_1=O(e^{-\lambda a_1})$.

In the case $|y|\leq b$, we will take $a=a_2>\sqrt{\pi/2}$
to ensure that the equation $Q(z,y)=0$ has
exactly two solutions $z_{\pm}=z_{\pm}(y)$ (with the same imaginary part) in the strip $0\leq \Im z\leq a_1$. 
The principal term of the asymptotics of the integral
\num{ifou} will come from the residue part of the RHS
in \num{intxres}.   

The three parameters $b$, $a_1$ and $a_2$ in the outlined
program are not determined rigidly and can be varied as long as we don't cross certain boundaries. Details are clarified in the following two lemmas. 

%%%%%%%%%%%%%%%%%%%%%%%%%%%%%%%%%%%%%%%
\iffalse
\begin{lem}\label{lem:ab}
For any $y\geq 0$, the system
$$
Q(x+iu,y)=0,\qquad 0\leq u< \sqrt{\f{3\pi}{2}}
$$ 
has exactly two real solutions $(x,u)$.
They have the same imaginary part and opposite real parts,
that is, the solutions are of the form $x=\pm x_*(y)$, $u=u_*(y)$.
The function $u_*(y)$ is continuous, increasing, and 
$$
 u_*(0)=\sqrt{\f{\pi}{2}}.
$$
\end{lem}

\pf
The complex quation $Q(x+iu,y)=0$ is equivalent to the pair
of real equations \num{zQ}.

Let $\gamma=\pi/2$ or $-\pi/2$. Substituting $x=\gamma/u$ into the equation $x^2+y^2=u^2$, we get
$$
 u^4-y^2 u^2-\gamma^2=0.
$$ 
Since $u$ is real, there must be $u^2\geq 0$, hence
$$
 u^2=\f{1}{2}\left(y^2+\sqrt{y^2+4\gamma^2}\right).
$$
Let $u_*(y)$ be the square root of the RHS.
Clearly, $u_*(y)$ is an increasing function.

If we replace the value $|\gamma|=\pi/2$ by $\pi(k+1/2)$
with $k\geq 1$, then $u^2\geq \gamma^2$ (by the above explicit formula). On the other hand,
$u^2<3\pi/2\leq|\gamma|$, a contradiction.
Hence the $|\gamma|=\pi/2$ is the only possibility
given the constraint on $u$. 

The rest is trivial. The value $\gamma=\pi/2$ yileds
the solution $u=u_*$, $x=x_*=\pi/(2u_*)$, while %the value 
$\gamma=-\pi/2$ yileds
the solution $u=u_*$, $x=-x_*$.
\eop
\fi
%%%%%%%%%%%%%%%%%%%%%%%%%%%%%%%%%%%%%%%

%%%%%%%%%%%%%%%%%%%%%%%%%%%%%%%%%%%%%%%
%\iffalse
\begin{lem}\label{lem:ab}
For any $a$ satisfying the inequalities
$$
 \sqrt{\f{\pi}{2}}<a<\sqrt{\f{3\pi}{2}},
$$
there exists a unique $b>0$ such that 
\\[1ex]
{\rm(i)} the equation
$Q(x+iu,y)=0$ does not have real solutions $(x,u,y)$ with
$0\leq u\leq a$ and $|y|>b$;
\\[1ex]
{\rm(ii)} for every $y\in(-b,b)$, the system
$$
Q(x+iu,y)=0,\qquad 0\leq u< a
$$ 
has exactly two real solutions $(x,u)$.
They have the same imaginary part and opposite real parts,
that is, the solutions are of the form $x=\pm x_*(y)$, $u=u_*(y)$.
\end{lem}

\pf 
We will prove that $b=\sqrt{a^2-x_0^2}$, where
$x_0=\pi/(2a)$. Note that the pair $x=x_0$, $y=b$ 
satisfies the system \num{zQ} with $u=a$, $k=0$.

\smallskip
(i) If $Q(x+iu,y)=0$ and $|y|>b$, $0\leq u\leq a$, then $x^2=u^2-y^2<a^2-b^2=x_0^2$.
Hence $|xu|<x_0 a=\pi/2$, and the equation $xu=\pi(2k+1)$
cannot hold.

\smallskip
(ii) %The situation is best seen graphically. 
Let $\gamma=\pi/2$ or $-\pi/2$. Substituting $x=\gamma/u$ into the equation $x^2+y^2=u^2$, we get
$$
 u^4-y^2 u^2-\gamma^2=0.
$$ 
Since $u$ is real, there must be $u^2\geq 0$, hence
\beq{explroot}
 u^2=\f{1}{2}\left(y^2+\sqrt{y^2+4\gamma^2}\right).
\eeq
Let $u_*(y)$ be the (positive) square root of the RHS.
Clearly, $u_*(y)$ is an increasing function.
By definition of the number $b$ (at the beginning of the proof), we have $u_*(b)=a$. Hence, $|y|\leq b$ implies
$0\leq u_*<a$, as required.

\smallskip
If we replace the value $|\gamma|=\pi/2$ by $\pi(k+1/2)$
with $k\geq 1$, then the three inequalities $u^2\geq \gamma^2$, $0\leq u\leq a$, and $ a^2<3\pi/2\leq |\gamma|$,
are incompatible. 

\smallskip
The rest is trivial. The value $\gamma=\pi/2$ yileds
the solution $u=u_*$, $x=x_*=\pi/(2u_*)$, while %the value 
$\gamma=-\pi/2$ yileds
the solution $u=u_*$, $x=-x_*$.
\eop
%\fi
%%%%%%%%%%%%%%%%%%%%%%%%%%%%%%%%%%%%%%%

\begin{lem}\label{aab}
There exist real numbers $b>0$ and $a_1$, $a_2$ satisfying the inequalities
\beq{a1a2}
 \sqrt{\f{\pi}{2}}<a_1<a_2<\sqrt{\f{3\pi}{2}},
\eeq
with the following properites: 
\\[1ex]
{\rm(i)} the equation
$Q(z,y)=0$ does not have solutions $(z,y)$ with
$y\in\RR$, $|y|\geq b$, 
$0\leq \Im z\leq a_1$;
\\[1ex]
{\rm(ii)} for every $y\in(-b,b)$, the system
$$
Q(z,y)=0, \qquad 0\leq \Im z<a_2
$$ 
has exactly two solutions $z_{\pm}(y)=\pm x_*(y)+i u_*(y)$.
\end{lem}
  
\pf
Choose $a_1$ and $a_2$ satisfynig the inequalities
\num{a1a2} arbitrarily. Take any $a\in(a_1,a_2)$
and determine $b$ as in Lemma~\ref{lem:ab}.
Let us check that the required properties are in place.

\smallskip
(i) If $|y|>b$, then by Lemma~\ref{lem:ab}, (i),
the equation $Q(z,y)$ does not have solutions 
with $0\leq \Im z <a$. 
The constraint $\Im z<a_1$ is even stronger. 
%In fact, with the latter constraint we
%can slightly relax the requirement $|y|>b$, changing
%it into $|y|>b_1$ with $b_1<b$. This covers the case $|y|=b$. 

\smallskip
(i) Let $|y|<b$. By Lemma~\ref{lem:ab}, (ii),  
the equation $Q(z,y)$ has two solutions $z=\pm x_*+iu_*$
satisfying the inequality $0\leq \Im z<a$. Relaxing
the constraint to  $0\leq \Im z<a_2$ could potentially
bring in extra solutions. 
For every such solution, $u=\Im z$ would be defined by \num{explroot} with $|\gamma|=\pi(k+1/2)$, $k\geq 1$.
Then $u^2\geq 3/2\pi$, which contradicts the condition
$u<a_2<\sqrt{3\pi/2}$.
\eop

%%%%%%%%%%
\begin{lem}
\label{lem:split}
Let $a_1$, $a_2$, $b$ be as in Lemma~\ref{aab} and
$I_1$, $I_2$ as in \num{isplit}.
Then 
\beq{iigta}
\ba{l}
\dst
I_1=O(e^{-\lambda a_1}),
\\[2ex]\dst
I_2=-\pi i\sum_{\pm}\int_{-b}^b \f{e^{i\lambda z_{\pm}(y)}}{z_{\pm}(y)}\,dy\;+\;O(e^{-\lambda a_2}),
\ea
\eeq
where $z_\pm(y)$ are the two solutions of the equation
$Q(z,y)=0$ defined in Lemma~\ref{aab}.
\end{lem}

\pf
1) Consider the integral \num{intt}
with $|y|>b$. By Lemma~\ref{aab} (i) and  \num{intxres} where we set $h=a_1$,
$$
T(y,\lambda)=\int_{-\infty}^\infty
\f{e^{\lambda (-a_1+ix)}}{Q(x+ia_1,y)}\,dx.
$$
Therefore,
$$
 |I_1|\leq \int_{|y|>b} |T(y,\lambda)|\,dy
 \leq e^{-\lambda a_1}\int_{|y|>b}\int_{-\infty}^\infty \f{dx\,dy}{|Q(x+ia_1,y)|}.
$$
The conditions of Lemma~\ref{lem:Qlb} are met with $u=a_1$,
$K=(-\infty,-b]\cup [b,\infty)$. 
The estimate \num{iigta} for $I_1$ follows at once.  

\smallskip
2) In \num{intxres} we set $h=a_2$. 
By Lemma~\ref{aab} (ii) there are two poles, $z=z_\pm(y)$, that contribute
to the sum of residues.
The residue at $z$ is
$$
 \f{e^{i\lambda z}}{Q'_z(z,y)}=
\f{e^{i\lambda z}}{2z e^{z^2+y^2}}=\f{e^{i\lambda z}}{-2z},
$$
since at a pole $Q(z,y)=0$ and $e^{z^2+y^2}=-1$.
We have identified the sum part of the RHS in \num{intxres}
with the first part in the RHS of the second estimate in 
\num{iigta}.

To show that the integral in the RHS of \num{intxres}
is $O(e^{-\lambda a_2})$, it 
suffices to repeat the same argument as in part 1 of this proof with obvious changes: $u=a_2$, $K=[-b,b]$.
\eop

\bigskip
To simplify the remaining calculations, note that
$$
 z_+=x_*+iu_*, \quad z_-=-x_*+iu_*,
$$
so $$iz_-=\overline{iz_+}.$$
Therefore
\beq{ccres}
-\pi i\sum_{\pm} \f{e^{i\lambda z_{\pm}}}{z_{\pm}}
\;=\;2\pi\,\Re\,\f{e^{i\lambda z_+}}{iz_+}.
\eeq

\begin{lem}
\beq{asres}
\int_{-b}^b \f{e^{i\lambda z_+(y)}}{iz_+(y)}
\;\sim\;-2^{1/2}\pi^{1/4}\,e^{-\lambda\sqrt{\pi/2}}\,
e^{i(\lambda\sqrt{\pi/2}+\pi/8)}
\eeq
as $\lambda\to+\infty$.
\end{lem}

\pf The function $\Re(i z_+(y))= -u_*(y)$ attains its
maximum $-\sqrt{\pi/2}$ at $y=0$, cf.~\num{explroot}.

The saddle point method tells us that
$$
\int_{-b}^b \f{e^{i\lambda z_+(y)}}{iz_+(y)}
\;\sim\;
\f{e^{i\lambda z_+(0)}}{iz_+(0)}\,\cdot\,\left(
\f{2\pi}{-i\lambda \,z_+''(0)}\right)^{1/2}.
$$ 
Recalling that $z_+^2=\pi i-y^2$, 
evaluate: 
$$
 iz_+(0)=i\sqrt{\pi} e^{i\pi/4}=\sqrt{\f{\pi}{2}}\,(-1+i)
$$
(since $\Re z_+, \Im z_+ >0$).

Expanding $z_+(y)$ in powers of $y^2$, we find  
$$
  z_+(y)=(\pi i)^{1/2}\left(1-\f{y^2}{2\pi i}\right)+O(y^4), 
$$
hence
$$
 z_+''(0)=\sqrt{\pi} e^{i\pi/4}\,\f{1}{-\pi i}
$$
and 
$$
-i z_+''(0)=\f{1}{\sqrt{\pi}} e^{i\pi/4}.
$$
Thus 
$$
\int_{-b}^b \f{e^{i\lambda z_+(y)}}{iz_+(y)}
\;\sim\;
\f{e^{(-1+i)\lambda\sqrt{\pi/2}}}{\sqrt{\pi} e^{3\pi i/4}}
\, \left(\f{2\pi^{3/2}}{\lambda e^{3\pi i/4}}\right)^{1/2}.
$$
Simplifying, we obtain \num{asres}.
\eop

\bigskip
Combining  \num{iigta}, \num{ccres}, \num{asres},
we get
$$
I_2\sim -2^{3/2}\pi^{5/4}\,e^{-\lambda\sqrt{\pi/2}}\;
\f{\cos(\lambda\sqrt{\pi/2}+\pi/8)}{\lambda^{1/2}},
$$
and, using \num{isplit},
\beq{anslam}
 S_*(\lambda)\;\sim\; 2^{3/2}\pi^{1/4}\,e^{-\lambda\sqrt{\pi/2}}\;
\f{\cos(\lambda\sqrt{\pi/2}+\pi/8)}{\lambda^{1/2}}.
\eeq
The final answer in the original notation, cf.~\num{tlam},
is
\beq{anst}
\mbox{\framebox{$\ba{c}\, \\[-1ex] \dst
S(t)\;\sim\; 
2\pi^{1/4}\,e^{-\sqrt{2\pi t}}\;\,
\f{\cos(\sqrt{2\pi t}+\pi/8)}{t^{1/4}}.
\\[-1ex] \,\ea
$}}
\eeq
\label{page:answer}

\bigskip
\noindent
\begin{picture}(200,200)
\put(0,0){\includegraphics[width=3in]{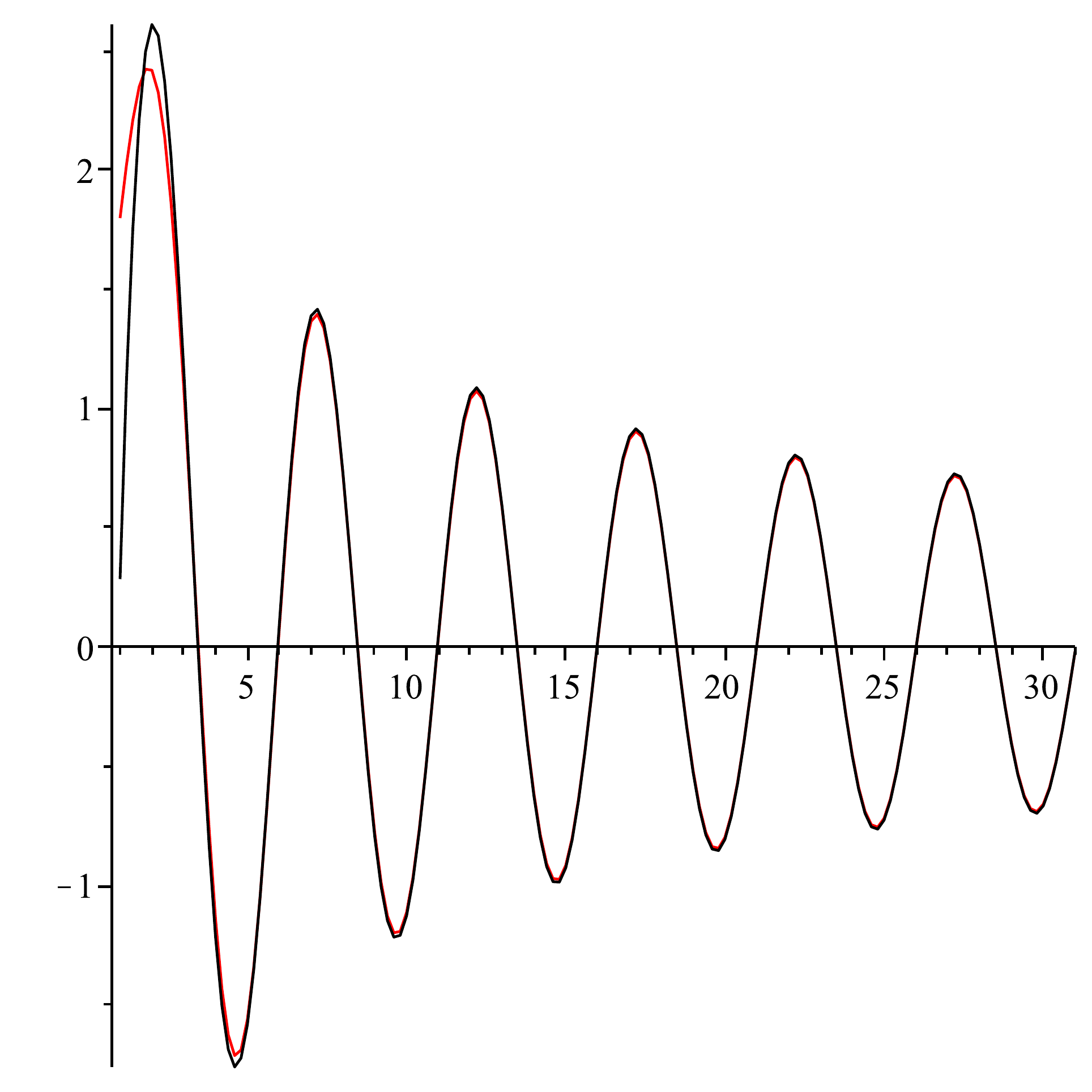}}
\put(210,178){\unl{Plots $\;-e^{\lambda\sqrt{\pi/2}} S_*(\lambda)\;$ vs $\;\lambda$}}
\put(220,148){Red: numerical quadrature \num{ireplam}}
\put(220,128){Black: asymptotics \num{anslam}}
\end{picture} 

\bigskip\bigskip
\noindent
{\bf Remark}.
From our analysis it is easy to see that the (unwritten) error term in 
\num{anslam} is
$$
 O\left(e^{-\lambda\sqrt{\pi/2}} \lambda^{-3/2}\right).
$$ 

The form in which the answer is presented in 
\num{anslam} or \num{anst}, 
is slightly inaccurate: the error term
needs not be dominated by the main term {\em everywhere},
since the latter becomes zero for some values of
$\lambda$ (or $t$).
 
The meaning of the $\sim$ sign is that the error term
is smaller by its {\it order of magnitude}, which is characterized by the {\it non-oscillating} (amplitude) factor
in  \num{anslam} or \num{anst}. 
Regarding the usage of the $\sim$ sign in this
and similar situations, I disagree with de Bruijn
who expressed negative opinion of such usage
({\it Asymptotic methods in analysis}, end of Sect.~5.11).

\end{document}